\documentclass[preprint,12pt]{elsarticle}

\usepackage{amssymb,amsthm}

\newtheorem{Def}{Definition}[section]
\newtheorem{Th}{Theorem}[section]
\newtheorem{Prop}{Proposition}[section]
\newtheorem{Not}{Notation}[section]
\newtheorem{Lemma}{Lemma}[section]

\newtheorem{Rem}{Remark}[section]
\newtheorem{Cor}{Corollary}[section]

\newtheorem{thm}{{\sc Theorem}}[section]

\newtheorem{lem}[thm]{{\sc Lemma}}

\newcommand\pa{\partial}

\newcommand\ep{\epsilon}

\newcommand\be{\begin{equation}}

\def\d{\begin{Def}}                     \def\t{\begin{Th}}
\def\p{\begin{Prop}}                    \def\n{\begin{Not}}
\def\la{\begin{Lemma}}                  \def\r{\begin{Rem}}
\def\c{\begin{Cor}}                     \def\ee{\begin{equation}}
\def\aa{\begin{eqnarray}}               \def\y{\begin{eqnarray*}}
\def\bd{\begin{description}}

\def\ed{\end{Def}}                      \def\et{\end{Th}}
\def\ep{\end{Prop}}                     \def\en{\end{Not}}
\def\el{\end{Lemma}}                    \def\er{\end{Rem}}
\def\ec{\end{Cor}}                      \def\eee{\end{equation}}
\def\eaa{\end{eqnarray}}                \def\ey{\end{eqnarray*}}
\def\ebd{\end{description}}

\def\nn{\nonumber}                      
\def\qe{\hfill{\rm Q.E.D.}}

                         \def\R{{\bf R}}

\journal{Journal of Mathematical Analysis and Applications}

\begin{document}

\begin{frontmatter}



\title{ Upper and lower bounds for normal derivatives  of spectral clusters of Dirichlet Laplacian }


\author{Xiangjin Xu \fnref{fn1}
}

\ead{xxu@math.binghamton.edu}

\address{Department of Mathematical Sciences, Binghamton University - SUNY\\ Binghamton, NY 13902, U.S.A}

\fntext[fn1]{The author is partially supported by NSF grants NSF-DMS-0602151 and NSF-DMS-0852507, and partially supported by Harpur College Grant in Support of Research, Scholarship and Creative Work in Year 2010-2011. }

\begin{abstract}
In this paper, we prove the upper and lower bounds for normal derivatives of spectral clusters $u=\chi_{\lambda}^s f$ of Dirichlet Laplacian $\Delta_M$,
$$c_s \lambda\|u\|_{L^2(M)} \leq \| \partial_{\nu}u \|_{L^2(\partial M)} \leq C_s \lambda \|u\|_{L^2(M)}  $$
where the upper bound is true for any Riemannian manifold, and the lower bound is true for some small $0<s<s_M$, where $s_M$ depends on the manifold only, provided that $M$ has no trapped geodesics (see Theorem \ref{Thm3} for a precise statement), which generalizes the early results  for single eigenfunctions by Hassell and Tao in 2002.


\end{abstract}


\end{frontmatter}

%\title{ Upper and lower bounds for normal derivatives  of spectral clusters of Dirichlet Laplacian}

%\author{ Xiangjin Xu\thanks{X. Xu is supported by the NSF Grant DMS-0602151 and DMS-0852507.}\\ Department of Mathematical Science, Binghamton University\\ Binghamton, NY 13902, U.S.A.\\ Email:xxu@math.binghamton.edu}
%\address{Department of Mathematical Science, Binghamton University, Binghamton, NY 13902}
%\email{xxu@math.binghamton.edu}
%\curraddr{ }

\date{}

%\maketitle

\section{\bf Introduction}

Let $M$ be a smooth compact Riemannian manifold with boundary $\partial M = Y$. It is well known that minus the Dirichlet Laplacian $-\Delta_M$ on $M$ has discrete spectrum $0 < \lambda_1^2 < \lambda_2^2 \leq \lambda_3^2 \dots \to \infty$. Let $e_j$ be an $L^2$-normalized eigenfunction corresponding to $\lambda_j^2$, and let $\psi_j$ be the normal derivative of $e_j$ at the boundary. In \cite{Ozawa} Ozawa posed the following question:
{\it Do there exist constants
$0< c < C < \infty$, depending on $M$ but not on $j$, such that
\begin{equation}
c \lambda_j \leq \| \psi_j \|_{L^2(Y)} \leq C \lambda_j  ?
\label{bounds}
\end{equation}}
Using heat kernel techniques, Ozawa \cite{Ozawa} showed that an averaged version of (\ref{bounds}) holds. More precisely, he showed that
$$
\sum_{\lambda_j < \lambda} \psi_j^2(y) =
\frac{\lambda^{n+2} }{ (4\pi)^{n/2} \Gamma((n/2)+2)}+ o(\lambda^{n+2}),
\quad \forall y \in Y.
$$
This asymptotic formula (after integrating over $Y$) would be implied by (\ref{bounds}) in view of Weyl asymptotics for the $\lambda_j$. In
\cite{HT} Hassell and Tao proved an upper bound of the form $\|\psi \|_2 \leq C\lambda$ for general manifolds, and a lower bound
$c \lambda \leq \|\psi \|_2$ provided that $M$ has no trapped geodesics ( see Theorem \ref{Thm3} for a precise statement).

%\vspace{2mm}

Define the  {\bf Spectral Cluster} $\chi_{\lambda}^s f$ with spectral band width $s$,
\begin{eqnarray*}
\chi_{\lambda}^s f &=&\sum_{\lambda_j\in[\lambda,\lambda+s)}e_j(f)
=\int_M[\sum_{\lambda_j\in[\lambda,\lambda+s)}e_j(x)e_j(y)]f(y)dy \label{SpecC} \\
e_j(f)&=& e_j(x)\int_M e_j(y)f(y)dy.
\end{eqnarray*}

The $L^2\to L^p$ estimates and gradient estimates on spectral clusters have been widely studied (see \cite{G}, \cite{SS1}-\cite{Xu2}). In general, the estimates for single eigenfunctions might be still true for spectral clusters in some sense. It is a natural question: {\it whether the upper and lower bound for normal derivatives of Dirichlet eigenfunctions in \cite{HT} is still true for normal derivatives of  spectral clusters $\chi_{\lambda}^s f$.}

The key obstacle to answer this question directly from the estimates of Hassell and Tao \cite{HT} for single eigenfunctions is that $\partial_{\nu}e_i$ and $\partial_{\nu}e_j$ are NOT orthogonal in $L^2(\partial M)$ in general when $\lambda_i\neq \lambda_j$.

%\vspace{5mm}

In this paper, based on the ideas in \cite{HT} plus some estimates for extra terms which come up for spectral clusters, we prove the upper bound from (\ref{bounds}) replacing $e_j$ by $\chi_{\lambda}^s f$ on general manifolds, for any $s>0$,

\t[{\bf Upper Bound}]\label{Thm1} Let $M$ be a smooth compact Riemannian manifold with boundary, and $u=\chi_{\lambda}^s f$ be the spectral clusters, we have that $\forall s>0$, there exists $C>0$ independent of $\lambda$ and $s$, such that
$$ \|\partial_{\nu} u\|_{L^2(\partial M)}\leq C\sqrt{1+s}(\lambda+s) \|u\|_{L^2(M)}.$$
\et

Especially for spectral projection $P_{\lambda}(f)=\displaystyle\sum_{\lambda_j\in(0,\lambda]}e_j(f)$, we have the upper bound estimate for its normal derivative:

\c\label{Cor1} Let $M$ be a smooth compact Riemannian manifold with boundary, for spectral projection $P_{\lambda}(f)$, we have
$$ \|\partial_{\nu}P_{\lambda}(f) \|_{L^2(\partial M)}\leq C\lambda^{3/2} \|P_{\lambda}(f)\|_{L^2(M)}.$$
\ec

Next it will be more subtle to study the lower bound from (\ref{bounds}) replacing $e_j$ by $\chi_{\lambda}^s f$. It might only hold when the value of $s$, the width of the spectral cluster, is sufficiently small. This can be seen in the case of the unit disc. If $s>\pi$, we can take a suitable linear combination of two consecutive eigenfunctions with angular dependence $e^{in\theta}$ (these are of the form $e^{in\theta} J_n(\alpha r)$ where $\alpha$ is a zero of the Bessel function $J_n$) and find a function in a ``wide" spectral cluster with zero normal derivative.

In order to obtain the lower bound, we first study on bounded Euclidean domains.  Following an idea of Rellich \cite{R} for single eigenfunctions on bounded Euclidean domains, we have the lower bound from (\ref{bounds}) replacing $e_j$ by $\chi_{\lambda}^s f$ on bounded Euclidean domains for small $s>0$,

\t[{\bf Lower Bound for Euclidean Domains}]\label{Thm2}
Let $M \subset \R^n$ be a bounded Euclidean domain, and $R_M=\max_{x,y\in M}|x-y|$ be the diameter of the domain $M$, and $u=\chi_{\lambda}^s f$ be the spectral clusters. Then for $0<s<\frac{1}{2R_M}$, there exists $C_s>0$ independent of $\lambda$, such that
$$ \| \partial_{ \nu} u\|_{L^2(\partial M)}\geq C_s\lambda \|u\|_{L^2(M)}.$$
\et

Next we turn to study the lower bound on general manifolds. To show a basic picture of our theorem, we refer some simple examples from \cite{HT} for single eigenfunctions, i.e., the cylinder (Example 3 in \cite{HT}), the hemisphere (Example 4 in \cite{HT}), the spherical cylinder (Example 5 in \cite{HT}). In all these examples, the upper bound holds, but the lower bound fails. These examples lead one to expect that the failure of the lower bound is related to the presence of geodesics in $M$ which do not reach the boundary. We obtain the lower bound estimates as in \cite{HT} replacing $e_j$ by $\chi_{\lambda}^s f$ for small $0<s<s_M$, where $s_M$ depends on the manifold only,

\t[{\bf Lower Bound for Manifolds}]\label{Thm3}
Suppose $M$ has {\bf no trapped geodesics}, i.e., $M$ can be embedded in the interior of a compact manifold with boundary, $N$, of the same dimension, such that every geodesic in $M$ eventually meets the boundary of $N$, and $u=\chi_{\lambda}^s f$ be the spectral clusters.   There exists $s_M>0$, which depends on the manifold only, such that for any $0<s<s_M$, there exists $C_s>0$ independent of $\lambda$, such that
$$ \|\partial_{\nu} u\|_{L^2(\partial M)}\geq C_s\lambda \|u\|_{L^2(M)}.$$
\et

We organize our paper as the following: In section 2, we prove a Rellich-type estimate from Green's formula, and some perturbation estimates to deal with the extra terms in the Rellich-type estimate. In section 3, we prove the the upper bound for general manifolds using the estimates from section 2 following the argument in \cite{HT}. In section 4, we prove the lower bound for Euclidean domains using the fact that the commutator $[-\Delta_M,x\cdot\nabla]=-2\Delta_M$, which gives the idea of the proof of lower bound for general case. In section 5, we show the lower bound for $L^2$ norm of $\partial_{\nu}u$ on an arbitrary compact Riemannian manifold $M$ satisfying the no trapped geodesics condition in Theorem \ref{Thm3} by finding a differential operator $P$ of order $2K-1$ which has a positive commutator with $-\Delta_M$, which depends on a trick due to Morawetz, Ralston and Strauss \cite{MRS}. In Appendix, we study the $L^2$ estimates for spectral clusters near the boundary, which are needed in the proof of Theorem \ref{Thm3}, following the same ideas in section 3 in \cite{HT} for single eigenfunctions.

In what follows we shall use the convention that $C$ denotes a constant that is not necessarily the same at each occurrence.

\section{\bf Rellich-type estimates and Perturbation estimates}

To prove the upper bound, and the lower bound for Euclidean domains, we use the following Lemma which we call a Rellich-type estimate.

\la({\bf Rellich-type estimates}) Let $u=\chi_{\lambda}^s f$ be the spectral projection of $f$. Then for any differential operator A,
\aa
\int_Y \partial_{\nu}  u Au d\sigma &= &\int_M <u,[-\Delta,A]u>dg + \int_M <(-\Delta-\lambda^2)u,Au>dg \nn\\
&&-\int_M <u,A(-\Delta-\lambda^2)u>dg.\label{Rellich}
\eaa
\el

{\bf Proof:} The proof is very simple. By Green's Formula, one has
\begin{eqnarray*} \int_Y \partial_{\nu}  u Au d\sigma - \int_Y u \partial_{\nu}  Au d\sigma=  \int_M <-\Delta u,Au>dg -\int_M <u,-\Delta Au>dg.\end{eqnarray*}
Note that $u\equiv 0$ on $Y$, left side of above equality gives left side of (\ref{Rellich}). Use the fact that
 $[-\Delta,A]=[-\Delta-\lambda^2,A]$ to write the right side as
\begin{eqnarray*}
&& \int_M <(-\Delta-\lambda^2)u,Au>dg -\int_M <u,(-\Delta-\lambda^2)Au>dg\\
&=& \int_M <u,[-\Delta,A]u>dg + \int_M <(-\Delta-\lambda^2)u,Au>dg \\
&&-\int_M <u,A(-\Delta-\lambda^2)u>dg.
\end{eqnarray*}\qe

\r If we pick $f=e_j$ the eigenfunction with eigenvalue $\lambda_j^2$, using the fact that $\Delta e_j+\lambda_j^2 e_j=0$, the above Lemma is reduced to Lemma 2.1 in \cite{HT}.
\er

Since there have two additional terms in (\ref{Rellich}) comparing with Lemma 2.1 in \cite{HT}, we need estimates them by the following Lemma.

\la({\bf Perturbation estimates}) Let $u=\chi_{\lambda}^s f$ be the spectral projection of $f$, $A$ is a differential operator with order one, we have
\begin{eqnarray*}
&&||(-\Delta-\lambda^2)u ||_2\leq 2s(\lambda+s) ||u||_2;\\
&& ||A(-\Delta-\lambda^2)u ||_2\leq C_A s(\lambda+s)^2 ||u||_2.
\end{eqnarray*}
\el

{\bf Proof:} For the first inequality, by direct computation, we have:
\begin{eqnarray*}
||(-\Delta-\lambda^2)u ||_2^2&=&\int_M<(-\Delta-\lambda^2)u ,(-\Delta-\lambda^2)u >dg\\
&=&\sum_{\lambda_j\in[\lambda,\lambda+s)}(\lambda_j^2-\lambda^2)^2e_j^2(f)\\
&<&\sum_{\lambda_j\in[\lambda,\lambda+s)}(2s\lambda+s^2)^2e_j^2(f)\\
&<&4s^2(\lambda+s)^2||u||_2^2
\end{eqnarray*}
For the second inequality, since $A$ is a differential operator with order one and $M$ is
compact, we have pointwise estimates
$$|A f(x)|\leq C_A|\nabla f(x)|,\qquad \forall x\in M\;  and \; \forall f\in C^1(M).$$
With this estimates, by direct computation, we have:
\begin{eqnarray*}
||A(-\Delta-\lambda^2)u ||_2^2&\leq& C_A^2||\nabla (-\Delta-\lambda^2)u ||_2^2\\
&=&C_A^2\int_M<\nabla(-\Delta-\lambda^2)u ,\nabla(-\Delta-\lambda^2)u >dg\\
&=&C_A^2\sum_{\lambda_j\in[\lambda,\lambda+s)}(\lambda_j^2-\lambda^2)^2\lambda_j^2e_j^2(f)\\
&<&C_A^2\sum_{\lambda_j\in[\lambda,\lambda+s)}(2s\lambda+s^2)^2\lambda_j^2e_j^2(f)\\
&<&C_A^24s^2(\lambda+s)^4||u||_2^2
\end{eqnarray*}
\qe

\section{\bf Upper bound for general manifolds}

In this section, we shall prove the the upper bound for general manifolds. Here we use the geodesic coordinates with respect to the boundary. We can find a small constant $\delta>0$ so that the map $x=(y,r)\in Y\times [0,\delta) \rightarrow M$,
sending $(y,r)$ to the endpoint $x$, of the geodesic of length $r$ which starts a $y\in Y=\partial M$ and is perpendicular to $Y$ is a local diffeomorphism. In this local coordinates $x=(y, r)$, the metric $ g = dr^2 + h_{ij} dy_i dy_j$ and the Riemannian measure
\begin{equation}
dg = k^2 dr dy, \quad where \quad k^4 = \det h_{ij}.\label{k}
\end{equation}
and the Laplacian can be written as
\begin{eqnarray*}
\Delta_g= \sum_{i,j=1}^n g^{ij}(x)\frac{\partial^2}{\partial x_i \partial x_j} + \sum_{i=1}^n b_i(x) \frac{\partial}{\partial x_i},
\end{eqnarray*}
where $(g^{ij}(x))_{1\le i,j \le n}$ is the inverse matrix of $(g_{ij}(x))_{1\le i,j \le n}$, and $g^{nn}=1$, and $g^{nk}=g^{kn}=0$ for $k\neq n$. Also the $b_i(x)$ are $C^{\infty}$ and real valued.

{\bf Proof of Theorem \ref{Thm1}:} For $u=\chi_{\lambda}^s f$, to prove an upper bound for the $L^2$ norm of $\partial_{\nu} u$, we choose an operator $A$ so that the left hand side of (\ref{Rellich}) in Lemma 2.1 is a positive form in $\partial_{\nu} u$. To do this, we choose $A = \chi(r) \partial_r$, where $\chi \in C_c^\infty(\R)$ is identically $1$ for $r$ close to zero, and vanishes for $r \geq \delta$. The left hand side of (\ref{Rellich}) in Lemma 3.1 is then precisely the square of the $L^2$  norm of $\partial_{\nu} u$.

After one integration by parts for the first term of the right hand side of (\ref{Rellich}) in Lemma 3.1, there are first order (vector-valued) differential operators $B_1$, $B_2$ with smooth coefficients,
$$\int_M <u,[-\Delta,A]u>dg=\int_M <B_1 u, B_2 u>dg. $$
From Lemma 3.2, each term of the right hand side of (\ref{Rellich}) in Lemma 3.1 is dominated by
\begin{eqnarray*}
|\int_M <u,[-\Delta,A]u>dg|&=&|\int_M<B_1 u, B_2 u>dg|\leq C_A||\nabla u||_2^2 \leq C_A(\lambda+s)^2||u||_2^2\\
|\int_M <(-\Delta-\lambda^2)u,Au>dg|&\leq& ||(-\Delta-\lambda^2)u||_2||Au||_2\leq C_As(\lambda+s)^2||u||_2^2\\
|\int_M <u,A(-\Delta-\lambda^2)u>dg|&\leq& ||A(-\Delta-\lambda^2)u||_2||u||_2\leq C_As(\lambda+s)^2||u||_2^2
\end{eqnarray*}
where $C$ and $C_A$ depend on the domain, but not on $\lambda$. This proves the upper bound for any compact Riemannian manifold with boundary.
\qe

If we choose $A = Q^*Q\partial_r$ near the boundary in the above proof, with $Q$ an elliptic differential operator of order k in the y variables, one has the $H^k$ estimates for upper bound of the spectral clusters $u=\chi_{\lambda}^s f$:
\t\label{Thm-H-k}
\aa
||\partial_{\nu} u||_{H^k(Y)}\leq C\sqrt{1+s}(\lambda+s)^k||u||_2\label{deriv-upper-bounds}
\eaa
for any integer k, and hence (by interpolation) any real k.
\et

\section{\bf Lower bound for Euclidean domains}

In this section, we shall prove Theorem \ref{Thm2}, the lower bound for Euclidean domain $M \subset \R^n$.

{\bf Proof of Theorem \ref{Thm2}:}
We choose $A$ so that the first term in the right hand side, rather than the left hand side, of (\ref{Rellich}) in Lemma 2.1 is a positive form. Without lose of generality, assume $M\subset \{x\in \R^n| |x|\leq \frac{R_M}{2}\}$.  We choose
\begin{equation}
A = \sum_{i=1}^n x_i \frac{\partial}{\partial x_i}=x\cdot\nabla.\nn
\end{equation}
As is very well known in scattering theory, the commutator of this with $-\Delta$ (which is minus the Euclidean Laplacian here) is $[-\Delta,A] = -2\Delta$, and for any $g\in C^1(M)$, $|A g(x)|\leq \frac{R_M}{2}|\nabla g(x)|$ for all $x\in M$. Hence, in this case the left side of (\ref{Rellich}) gives us
\begin{equation}
\int\limits_{Y} \frac{\partial u}{\partial \nu} Au \, d\sigma =
\int\limits_{Y} \nu \cdot x \, \big( \frac{\partial u}{\partial \nu} \big)^2 \, d\sigma
\leq C \| \partial_{\nu}u \|_2^2,
\label{lower-Euc}
\end{equation}
And the right side of  (\ref{Rellich}) gives us
\begin{eqnarray}
&&RIGHTSIDE \; of\;  (\ref{Rellich}) \nn\\
&\geq& \int_M <u,-2\Delta u>dg -\|(-\Delta-\lambda^2)u\|_2\|Au \|_2-\|u\|_2\|A(-\Delta-\lambda^2)u\|_2\nn\\
&\geq& 2\|\nabla u\|_2^2-\frac{R_M}{2}\Big[\|(-\Delta-\lambda^2)u\|_2\|\nabla u\|_2+\|u\|_2\|\nabla(-\Delta-\lambda^2)u\|_2\Big] \nn\\
&\geq& \Big(2\lambda^2-2R_M s(\lambda+s)^2\Big) \|u\|_2^2,\label{lower-Euc-2}
\end{eqnarray}
which gives the lower bound. The equality in (\ref{lower-Euc}) for single eigenfunctions was proved by Rellich \cite{R}.
\qe

\section{\bf The lower bound on Riemannian manifolds}

To find a lower bound for $L^2$ norm of $\partial_{\nu}\Big(\chi_{\lambda}^s f\Big)$ on an arbitrary compact Riemannian manifold $M$ satisfying the no trapped geodesics conditions of the main Theorem, we need to find a differential operator which has a positive commutator with $-\Delta_M$ as we did for domains in Euclidean spaces. One might wonder whether, on an arbitrary compact Riemannian manifold, with no trapped geodesics, one could choose a first order {\it differential} operator $A$ whose commutator with $-\Delta_M$ had a positive symbol.  Example 8 in \cite{HT} shows  that this is impossible in general.

Firstly, we have a first order pseudo-differential operator $A$ on $N$ which has the required property to leading order, i.e., such that the symbol of $i[-\Delta,A]$ is positive:

\begin{lem}\label{Q-lemma}({\bf Lemma 4.1 in \cite{HT}}) Given any geodesic $\gamma$ in $S^*N$, there is a first order, classical, self-adjoint pseudodifferential operator $Q$ satisfying the transmission condition (see \cite{Ho}, section 18.2), and properly supported on $N$, such that the principal symbol $\sigma(i[-\Delta,Q])$ of $i[-\Delta,Q]$ is nonnegative on $T^*M$, and
\begin{equation}
\sigma(i[-\Delta,Q]) \geq \sigma(-\Delta) = |\xi|^2
\label{comm-cond}
\end{equation}
on  a  conic  neighborhood  $U_\gamma$ of  $\gamma \cap T^*M$.
\end{lem}

We now use Lemma~\ref{Q-lemma} to construct our operator $A$. For each geodesic $\gamma$ in $S^*N$, we have a conic neighborhood $U_\gamma$ as in the Lemma. By compactness of $S^*M$, a finite number of the $U_\gamma$ cover $S^*M$. Let $A$ be the sum of the corresponding $Q_\gamma$. Then Lemma~\ref{Q-lemma} implies that
\begin{equation}
\sigma(i[-\Delta,A]) \geq |\xi|^2 \; on \; T^*M .
\label{A}
\end{equation}

Secondly, we turn the pseudodifferential operator $A$ into a differential operator $P$  of order $2K-1$ with positive commutator with $-\Delta$  as  Hassell and Tao did at Section 5 in \cite{HT} for single eigenfunctions, which depends on a trick due to Morawetz, Ralston and Strauss \cite{MRS}.

Recall some facts about spherical harmonics. Let $\Delta_{S^{n-1}}$ denote the Laplacian on the $(n-1)$-sphere, which has eigenvalues $k(n+k-2)$, $k = 0, 1, 2, \dots$, and the corresponding eigenspace be denoted $V_k$. We recall that  for every $\phi \in V_k$, the function $r^k \phi$ (thought of as a function on $\R^n$ written in polar coordinates) is a homogeneous polynomial, of degree $k$, on $\R^n$. We summarize the needed results from Section 5 in \cite{HT} as the following proposition:

\p[Hassell-Tao \cite{HT}]
Since the symbol $a$ of the operator $A$  is odd, there is spherical harmonics expansion of $a$ restricted to the cosphere bundle of $N$
$$
a \restriction_{ S^*N}= \sum_{l=0}^\infty  \phi_{2l+1}(x, \frac{\xi}{|\xi|}), \quad \phi_k(x, \cdot) \in V_k(S^*_x N).
$$
And  there is a nature number $K$ such that the operator $A'$ with symbol
$$
a' = \sum_{l=0}^{K-1} \phi_{2l+1}(x, \frac{\xi}{|\xi|})
$$
also has positive commutator with $-\Delta$. Following \cite{MRS}, one can turn $A'$ into a differential operator $P$ of order $2K-1$, by letting
$$
p = \sigma(P) = \sum_{l=0}^{K-1} \phi_{2l+1}(x, \frac{\xi}{|\xi|}) |\xi|^{2K- 1}.
$$
Moreover, the symbol of $i[-\Delta,P]$ satisfies
$$
\sigma(i[-\Delta,P]) = |\xi|^{2K} \big( \sigma(i[-\Delta,P]) |_{|\xi| = 1} \big) \geq c|\xi|^{2K}
\quad  for\; some\; c > 0.
$$
Applying the G$\mathring{a}$rding inequality to $Q = i[-\Delta,P]$, there is
\ee
\int_M \langle u, Qu \rangle dg \geq c \| u \|_{H^K(M)}^2 - C \Big(  \| u \|_{L^2(M)}^2 + \sum_{k=0}^{K-1} \| \pa_r^k u \| _{H^{K-1/2-k}(Y)}^2 \Big),
\label{Garding}
\eee
where $c$ is a positive constant depending on $P$ and $(M,g)$.
\ep

{\bf Proof of Theorem \ref{Thm3}:}
Let $A=P$ in Lemma 2.1,
\aa
&&RIGHTSIDE \; of\;  (\ref{Rellich}) \nn\\
&\geq& \int_M <u,Qu>dg -\|(-\Delta-\lambda^2)u\|_2\|Pu \|_2-\|u\|_2\|P(-\Delta-\lambda^2)u\|_2\nn\\
&\geq&  \int_M <u,Qu>dg -C\Big[\|(-\Delta-\lambda^2)u\|_2\|u\|_{H^{2K-1}(M)}\nn\\
&& +\|u\|_2\|(-\Delta-\lambda^2)u\|_{H^{2K-1}(M)}\Big] \nn\\
&\geq&  \int_M <u,Qu>dg -Cs\lambda^{2K}\|u\|_2^2\nn\\
&\geq&  (c-Cs)\lambda^{2K} \| u \|_2^2 - C \Big(  \| u \|_{L^2(M)}^2 + \sum_{k=0}^{K-1} \| \pa_r^k u \| _{H^{K-1/2-k}(Y)}^2 \Big), \label{lower-general}
\eaa
where we use Lemma 2.2 and $\|u\|_{H^k(M)}\leq C\lambda^k\|u\|_2,\; \forall k>0$ to estimate the extra terms, and make use of (\ref{Garding}) to obtain the last inequality. Thus, there exists a constant $s_M>0$, which depends on $M$ only, such that the first term in (\ref{lower-general}) is positive when $0<s<s_M$.

Next to consider the left hand side of (\ref{Rellich}). Let us write $u = k^{-1}v$, where $k$ is as in (\ref{k}), so that $v$ satisfies (\ref{v-eqn}) in Appendix. Then $Pu = \tilde P v$, where $\tilde P = P \circ k$ is a differential operator of order $2K-1$. Since $k$ is smooth, we obtain
\begin{equation}
C_s\lambda^{2K}\|u\|_2^2 \leq \|u\|_2^2 + \sum_{k=0}^{K-1} \| \pa_r^k v \| _{H^{K-1/2-k}(Y)}^2
+  \big| \int_{Y} \langle \partial_{\nu}v, \tilde P v \rangle \, d\sigma \big|.
\label{eqqq}\end{equation}
Since $v=0$ at $Y$, and we are interested in $\tilde P v |_Y$, we may assume that $\tilde P = P' \pa_r$, where $P'$ has order $2K-2$. Using (\ref{v-eqn}) in Appendix, we may replace $\pa_r^2 v$ by $-(\lambda - F)v - \pa_{y_i} (h^{ij} \pa_{y_j} v)+H$ repeatedly, until only $\pa_r \pa_y^\alpha v$ terms remain. Thus we have
\begin{eqnarray*}
\tilde P v |_Y = \sum_{j=0}^{K-1} \lambda^j P_j (\partial_{\nu}v),
\end{eqnarray*}
where $P_j$ is a differential operator on $Y$ of order $2(K-1-j)$, independent of $\lambda$. Hence (\ref{eqqq}) becomes
\begin{equation}
C'_s\lambda^{2K}\|u\|_2^2 \leq \|u\|_2^2 + \sum_{k=0}^{K-1} \lambda^{2k-2} \| \partial_{\nu}u \| _{H^{K-1/2-k}(Y)}^2
+  \sum_{j=0}^{K-1} \lambda^{2j}\big| \int_{Y} \langle \partial_{\nu}v, P_j (\partial_{\nu}v) \rangle \, d\sigma \big|.\nn
\end{equation}
The argument to reduce $Pu$ on $Y$ to $\sum_{j=0}^{K-1} \lambda^j P_j (\partial_{\nu}v)$ on $Y$ is the same as what Hassell and Tao did at Section 5 in \cite{HT} for single eigenfunctions.

Using the upper bound estimate (\ref{deriv-upper-bounds}) for $H^k$ norm on the sum over $k$ and for all terms in the sum over $j$ with $j < K-1$, we find
\begin{eqnarray*}
C''_s\lambda^{2K}\|u\|_2^2 \leq  (1 + \lambda^{2K-1})\|u\|_2^2 + \lambda^{2K-2} \|  \partial_{\nu}u \|_2^2
+  \lambda^{2K-1} \|u\|_2\|  \partial_{\nu}u\|_2 .
\end{eqnarray*}
which gives
\begin{equation}
\| \partial_{\nu}u \|_2^2+\lambda\|u\|_2\|  \partial_{\nu}u\|_2-(C''_s-\lambda^{-1}-\lambda^{-2K} )\lambda^{2}\|u\|_2^2\geq 0.\label{almost}
\end{equation}
Solve the inequality (\ref{almost}), for $\lambda$ large enough, we have constant $C_s$ independent of $\lambda$, such that
\begin{eqnarray*}
\| \partial_{\nu}u \|_2\geq C_s\lambda\|u\|_2,
\end{eqnarray*}
This proves the lower bound.
\qe

%\vspace{2mm}

\r
One may also prove the lower bound following what Hassell and Tao did at Section 4 in \cite{HT} for single eigenfunctions almost line by line, while one need do some additional estimates on the nonhomogeneous terms like $H$ in (\ref{v-eqn}), which can be looked as small perturbation terms when $s>0$ is small enough. This approach is length and involves many pseudodifferential operator constructions and calculus.
\er

\section{\bf Appendix: Estimates for spectral clusters near the boundary}

Here we study the $L^2$ estimates for spectral clusters near the boundary, which are needed in the proof of Theorem \ref{Thm3}, following some ideas from section 3 in \cite{HT} with its erratum \cite{HT1} for single eigenfunctions and the upper bound for $\partial_{\nu}\Big(\chi_{\lambda}^s f\Big)$ from Theorem \ref{Thm1}.

As in section 3, we use the geodesic coordinates system $(y,r)$ near the boundary. Let us denote the boundary of $M$ by $Y$, and write $Y_r$ for the set of points at distance $r$ from the boundary, which is a submanifold for $r \le \delta$. Suppose that $u=\chi_{\lambda}^s f$ is a spectral cluster for Dirichlet Laplacian. Similar as Lemma 3.2 in \cite{HT} with its erratum \cite{HT1} for a single eigenfunction, we derive an estimate on the $L^2$ norm of the spectral cluster $u=\chi_{\lambda}^s f$ on $Y_r$, exploiting the fact that $u$ vanishes on the boundary.

\begin{Prop} \label{L7-2}  There exists $C > 0$, independent of $\lambda$ and $s$, such that
\begin{equation}
\int_{Y_r} u^2 d\sigma(y) \leq C\sqrt{1+s}(\lambda+s)^2 r^2\|u\|_2^2 \quad \forall\; r \in [0, \frac{\delta}{3}].
\label{bdy-est}\end{equation}
\end{Prop}

It will be convenient to change to the function $v = ku$ (this is equivalent to looking at the Laplacian acting on half-densities). Denote $u_l=e_l(f)$, $v_l=ku_l$ for $l\in[\lambda, \lambda+s)$. From the equation (3.2) in \cite{HT}, $v_l$ solves the equation
\begin{equation}
\pa_r^2 v_l + \pa_i (h^{ij} \pa_j v_l) + \lambda_l^2 v_l + F v_l = 0, \quad  h^{ij} = (h_{ij})^{-1}\nn
\end{equation}
where
$$
F = - k^{-1}\pa_r^2 k - k^{-1} \pa_i ( h^{ij} \pa_j k )
$$
is a smooth function on $M$. We have
\begin{equation}\label{v-eqn}
\pa_r^2 v + \pa_i (h^{ij} \pa_j v) + \lambda^2 v + F v = \sum_{\lambda_l\in[\lambda, \lambda+s)}(\lambda^2-\lambda_l^2)v_l=H
\end{equation}

As did in section 3 of \cite{HT} for a single eigenfunction, where the nonhomogeneous term $H$ doesn't appear, for the spectral cluster $u$, we define a sort of `energy' $E(r)$ for each value of $r$:
\begin{equation}
E(r) = \frac1{2} \int_{Y_{r}} \big(  v_r^2 + (\lambda^2 + F)v^2 - h^{ij} \pa_i v \pa_j v-Hv \big) dy.
\end{equation}
This is obtained formally from the energy for hyperbolic operators, with $r$ playing the role of a time variable, by switching the sign of the term involving tangential derivatives. Similar as Lemma 3.1 in \cite{HT}, we have the following estimate for $E(r)$

\begin{lem}\label{L7-1} For $r \in [0, \delta]$,
\begin{equation}\label{energy-est}
|E(r)| \leq C\sqrt{1+s}(\lambda+s)^2\|u\|_2^2
\end{equation}
where $C$ is independent of $\lambda$ and $s$.
\end{lem}

{\bf Proof of Lemma \ref{L7-1}:}
From the upper bound argument in section 3, we know that $E(0) = \frac1{2} \| \partial_{\nu} u \|_2^2 \leq C\lambda^2\| u \|_2^2$. We compute the derivative of $E(r)$:
\begin{eqnarray*}
\frac{\pa}{\pa r} E(r) = \int_{Y_{r}} \big( \pa_r^2 v \pa_r v +(\lambda^2 + F)v \pa_r v  -h^{ij} \pa_i v \pa_r \pa_j v -Hv_r\\
+ \frac{\pa h^{ij}}{\pa r} \pa_i v \pa_j v + \frac{\pa F}{\pa r} v^2 -H_rv\big) dy.
\end{eqnarray*}
Integrating by parts in the third term, using the equation for $v$, and applying Cauchy-Schwarz to last term, we
obtain
\begin{eqnarray*}
\big| \frac{\pa E}{\pa r} \big|(r) &\leq& C \int_{Y_{r}} \big( v^2 + |\nabla v|^2 +\lambda^2|v|^2+\lambda^{-2}|\nabla H|^2\big) dy \\
&\leq& C \int_{Y_{r}} \big( u^2 + |\nabla u|^2+ \lambda^2|u|^2+\lambda^{-2}|\nabla (H/k)|^2\big) k^2 dy.
\end{eqnarray*}
Thus, for $r_0 \in [0, \delta]$,
\begin{eqnarray*}
E(r_0) &=& E(0) + \int_0^{r_0} \frac{d}{dr} E(r) dr  \\
&\leq & C\lambda^2\| u \|_2^2 + \int_M \big( u^2 + |\nabla u|^2+ \lambda^2|u|^2
+\lambda^{-2}|\nabla (H/k)|^2\big) dg\\
&\leq&  C\sqrt{1+s}(\lambda+s)^2\|u\|_2^2,
\end{eqnarray*}
where we use Lemma 2.2 to estimate them last term.
\qe

\vspace{4mm}

{\bf Proof of Proposition \ref{L7-2}:}
Here we follow the main idea of proof of Lemma 3.2 in \cite{HT} with its erratum \cite{HT1} for single eigenfunctions.

Consider the $L^2$ norm on $Y_r$,
$$L(r) = \int_{Y_r} u^2 k^2 dy = \int_{Y_r} v^2 dy.$$
And we have
\begin{equation}\label{u-bd}
\int_0^\delta L(r)  dr \leq  \int_M u^2 \, dg = \|u\|_2^2.
\end{equation}
By direct computation, we have
$$ L'(r_0) = \int_{Y_{r_0}} 2 v v_r\ dy\quad \mathrm{and} \quad  L''(r_0) =4 \int_{Y_{r_0}}  v_r^2  \ dy \ -  4E(r_0).$$
On the other hand, from Cauchy-Schwarz we have
$$ 4 \int_{Y_{r_0}} v_r^2 \ \geq \ \frac{(\int_{Y_{r_0}} 2 v v_r\ dy)^2}{\int_{r=r_0} v^2\ dy} = \frac{L'(r_0)^2}{L(r_0)}.$$
Thus we have the differential inequality for $L(r)$:
\begin{equation}
L'' \geq \frac{(L')^2}{L} - C\sqrt{1+s}(\lambda+s)^2\|u\|_2^2 ,
\label{diff-ineq}
\end{equation}
for some constant $C$ depending only on the manifold $M$. Define the quantity
$$
B(r) := \frac{L'(r)^2}{L(r)^2} - \frac{C\sqrt{1+s}(\lambda+s)^2\|u\|_2^2}{L(r)}
$$
For any $r\in [0, \delta]$ with $L'(r)>0$, from (\ref{diff-ineq}) we have
$$
B'(r) = \frac{2 L' L''}{L^2} - \frac{2 (L')^3}{L^3} + \frac{2 C\sqrt{1+s}(\lambda+s)^2\|u\|_2^2 L'}{L^2} \geq 0.
$$
Hence $B(r)$ is non-decreasing for $r\in [0, \delta]$ with $L'(r)>0$.

\vspace{4mm}

{\bf Claim:} There is $\Lambda>0$ such that for any $\lambda\geq \Lambda$, either $L'(r)\leq 0$ or $B(r)\leq 0$ are true for all $0<r<\delta/3$.

\vspace{4mm}

Define $O_{\lambda}=\{r | L'(r)>0,\; 0<r<\delta \}=\cup (a_n^{\lambda}, b_n^{\lambda})\subset (0,\delta]$. For any $r\in (a_n^{\lambda}, b_n^{\lambda})$ with $b_n^{\lambda}<\delta$, we have 
$$
B(r)\leq B(b_n^{\lambda})= - \frac{C\sqrt{1+s}(\lambda+s)^2\|u\|_2^2}{L(b_n^{\lambda})}\leq 0.
$$

Hence {\bf Claim} will be true unless there is an unbounded sequence of $\lambda$ such that $B(r_0)> 0$ for some $r_0\in (a_n^{\lambda}, \delta)$ and $0 < r_0 < \delta/3$, where $(a_n^{\lambda}, \delta)$ is one subinterval of $O_{\lambda}$.  Then we would have $B(r) > 0$ for all $r \geq r_0$, so
$$ 
L'(r)^2 > 2C\sqrt{1+s}(\lambda+s)^2\|u\|_2^2 L(r) \hbox{ for all } r \geq r_0.
$$
In particular $L'(r)$ would be strictly positive for $r \geq r_0$. We rearrange this as
$$ 
(L(r)^{1/2})' = \frac{1}{2} L'(r) L(r)^{-1/2} > \sqrt{C\sqrt{1+s}(\lambda+s)^2\|u\|_2^2} \hbox{ for all } r > r_0.
$$
This would give
$$ 
L(r) \geq C'(1+s) \lambda^2(r-r_0)^2\|u\|_2^2 \hbox{ for all } r >\delta/3>r_0.
$$
This would contradict the bound (\ref{u-bd}) for large $\lambda$. Hence {\bf Claim} is true.

\vspace{4mm}

Thus, for $\lambda \geq \Lambda$ (where $\Lambda$ is obtained from {\bf Claim} and depends only on $M$), we must have $B(r) \leq 0$ or $L'(r) \leq 0$ for all $0 < r  \leq \delta/3$. In either case
$$ 
(L(r)^{1/2})' = \frac{1}{2} L'(r) L(r)^{-1/2} \leq \sqrt{C\sqrt{1+s}(\lambda+s)^2\|u\|_2^2} \hbox{ for all } r \leq \frac{\delta}{3}.
$$
Since $L(0) = 0$, this implies (\ref{bdy-est}) for $\lambda \geq \Lambda$.

Next for $\lambda < \Lambda$, since 
$$
u=\sum_{\lambda_j\in[\lambda,\lambda+s)}e_j(f)\leq \Big(\sum_{\lambda_j\in[\lambda,\lambda+s)}e_j^2\Big)^{1/2}\Big(\sum_{\lambda_j\in[\lambda,\lambda+s)}||e_j(f)||_2^2\Big)^{1/2}
=\Big(\sum_{\lambda_j\in[\lambda,\lambda+s)}e_j^2\Big)^{1/2}||u||_2
$$
we have
$$
\int_{Y_r} u^2 d\sigma(y)\leq \int_{Y_r} \sum_{\lambda_j\in[\lambda,\lambda+s)}e_j^2 d\sigma(y)\|u\|_2^2 \leq C\Big(\sum_{\lambda_j\in[\lambda,\lambda+s)}\lambda_j^2\Big) r^2\|u\|_2^2\leq C_{\Lambda} r^2\|u\|_2^2
$$
where we make use of the result of Lemma 3.2 in \cite{HT} with its erratum \cite{HT1} for single eigenfunctions:
$$
\int_{Y_r} e_j^2 d\sigma(y)\leq C\lambda_j^2, \quad \lambda_j\leq \Lambda+s.
$$
\qe

\vspace{5mm}

{\bf Acknowledgement:} The author would like to thank Professor Andrew Hassell for pointing out a mistake in first version of this paper and some helpful suggestions on this paper.

\bibliographystyle{abbrv}

\end{document}